\documentclass[11pt]{amsart}
\usepackage{amssymb}
\newtheorem{thm}{Theorem}
\newtheorem{cor}{Corollary}
\theoremstyle{definition}
\newtheorem*{dfn}{Definition}
\newtheorem{rem}{Remark}
\newcommand{\rhobin}{\mathbin{\rho}}
\newcommand{\barrhobin}{\mathbin{\bar\rho}}
\begin{document}
\title[Insertion of a contra-continuous function]{A sufficient condition 
for the insertion of a contra-continuous (Baire-one) function}
\author{Majid Mirmiran}
\address{Department of Mathematics,\\
University of Isfahan\\
Isfahan 81744, Iran.}
\email{mirmir@sci.ui.ac.ir}
\subjclass[2000]{Primary 54C08, 54C10; Secondary 26A15, 54C30}
\keywords{Contra-continuous function, Baire-one function,
$\Lambda$-sets, Lower cut set}
\begin{abstract}
A sufficient condition for the insertion of a contra-continuous
(resp.\ Baire-one) function between two comparable real-valued
functions is given on the topological spaces that $\Lambda$-sets
are open (resp.\ $G_{\delta}$-sets).
\end{abstract}
\maketitle

\section{Introduction}

Results of Kat\v{e}tov [4], [5] concerning binary relations and the
concept of an indefinite lower cut set for a real-valued function, which
is due to Brooks [1], are used in order to give a sufficient condition for
the insertion of a contra-continuous (resp.\ Baire-one) function between
two comparable real-valued functions on the topological spaces that
$\Lambda$-sets [7] are open (resp.\ $G_{\delta}$-sets).

A generalized class of closed sets was considered by Maki in 1986 [7]. He
investigated the sets that can be represented as union of closed sets and
called them $V$-sets. Complements of $V$-sets, i.e., sets that are
intersection of open sets are called $\Lambda$-sets [7].

A real-valued function $f$ defined on a topological space $X$ is called
{\it contra-continuous} [2] (resp.\ {\it Baire-one}) if the preimage of
every open subset of $\mathbb{R}$ is closed (resp.\ $F_{\sigma}$-set) in
$X$.

If $g$ and $f$ are real-valued functions defined on a space $X$, we write
$g\leq f$ in case $g(x)\leq f(x)$ for all $x$ in $X$.

\section{The main result}

Before giving a sufficient condition for insertability of a
contra-continuous (Baire-one) function, the necessary definitions and
terminology are stated.

\begin{dfn}
Let $A$ be a subset of a topological space $(X, \tau)$. 
We define the subsets $A^{\Lambda}$ and $A^V$ as follows: 
$A^{\Lambda} = \bigcap\{O: O\supseteq A, O\in (X, \tau)\}$ and
$A^{V} = \bigcup\{F:F\subseteq A, F^c\in (X, \tau)\}$.
In [3], [6], $A^{\Lambda}$ is called the {\it kernel} of $A$.
\end{dfn}

The following first two definitions are due to, or are modifications of, 
conditions considered in [4], [5].

\begin{dfn}
If $\rhobin$ is a binary relation in a set $S$ then $\barrhobin$ is 
defined as follows: $x \barrhobin y$ if and only if $y \rhobin \nu$ 
implies $x \rhobin \nu$ and $u \rhobin x$ implies $u \rhobin y$ for any 
$u$ and $v$ in $S$.
\end{dfn}

\begin{dfn}
A binary relation $\rhobin$ in the power set $\mathcal{P}(X)$  of a 
topological space $X$ is called a {\it strong binary relation} in 
$\mathcal{P}(X)$ in case $\rho$ satisfies each of the
following conditions:
\begin{enumerate}
\item
If $A_i \rhobin B_j$ for any $i\in \{1,\dots,m\}$ and for any
$j\in\{1,\dots,n\}$, then there exists a set $C$ in $\mathcal{P}(X)$ such
that $A_i \rhobin C$ and $C \rhobin B_j$ for any $i\in \{1,\dots,m\}$
and any $j\in \{1,\dots,n\}$.
\item
If $A\subseteq B$, then $A \barrhobin B$.
\item
If $A \rhobin B$, then $A^{\Lambda} \subseteq B$ and $A\subseteq B^V$.
\end{enumerate}
\end{dfn}

The concept of a lower indefinite cut set for a real-valued
function was defined by Brooks [1] as follows:

\begin{dfn}
If $f$ is a real-valued function defined on a space $X$ and if 
$\{x\in X: f(x) < l\} 
\subseteq A(f,l) 
\subseteq \{x\in X: f(x)\leq l\}$ 
for a real number $l$, then $A(f,l)$ is a {\it lower indefinite cut 
set} in the domain of $f$ at the level $l$.
\end{dfn}

We now give the following main result:

\begin{thm} %\noindent {\bf Theorem 1.} 
Let $g$ and $f$ be real-valued
functions on the topological space $X$, that $\Lambda$-sets in
$X$ are open (resp.\ $G_{\delta}$-sets), with $g\leq f$. If there
exists a strong binary relation $\rhobin$ on the power set of $X$
and if there exist lower indefinite cut sets $A(f,t)$ and
$A(g,t)$ in the domain of $f$ and $g$ at the level $t$ for each
rational number $t$ such that if $t_1<t_2$ then $A(f,t_1) \rhobin
A(g, t_2)$, then there exists  a contra-continuous (resp.
Baire-one) function $h$ defined on $X$ such that $g\leq h\leq f$.
\end{thm}

\begin{proof}
Let $g$ and $f$ be real-valued functions defined on
the $X$ such that $g\leq f$. By hypothesis there exists a strong
binary relation $\rhobin$ on the power set of $X$ and there exist
lower indefinite cut sets $A(f,t)$ and $A(g,t)$ in the domain of
$f$ and $g$ at the level $t$ for each rational number $t$ such
that if $t_1<t_2$ then $A(f,t_1) \rhobin A(g, t_2)$.

Define functions $F$ and $G$ mapping the rational numbers
$\mathbb{Q}$ into the power set of $X$ by $F(t)=A(f,t)$ and
$G(t)=A(g,t)$. If $t_1$ and $t_2$ are any elements of
$\mathbb{Q}$ with $t_1<t_2$, then $F(t_1) \barrhobin F(t_2)$,
$G(t_1) \barrhobin G(t_2)$, and $F(t_1) \rhobin G(t_2)$. By Lemmas 1
and 2 of [5] it follows that there exists a function $H$ mapping
$\mathbb{Q}$ into the power set of $X$ such that if $t_1$ and
$t_2$ are any rational numbers with $t_1<t_2$, then
$F(t_1) \rhobin H(t_2)$, $H(t_1) \rhobin H(t_2)$ and 
$H(t_1) \rhobin G(t_2)$.

For any $x$ in $X$, let $h(x)=\inf\{t\in \mathbb{Q}: x\in H(t)\}$.

We first verify that $g\leq h\leq f$: If $x$ is in $H(t)$ then $x$
is in $G(t')$ for any $t'>t$; since $x$ in $G(t')=A(g,t')$
implies that $g(x)\leq t'$, it follows that $g(x)\leq t$. Hence
$g\leq h$. If $x$ is not in $H(t)$, then $x$ is not in $F(t')$
for any $t'<t$; since $x$ is not in $F(t')=A(f,t')$ implies that
$f(x)>t'$, it follows that $f(x)\geq t$. Hence $h\leq f$.

Also, for any rational numbers $t_1$ and $t_2$ with $t_1<t_2$, we
have $h^{-1}(t_1, t_2)= H(t_2)^V\setminus H(t_1)^{\Lambda}$.
Hence $h^{-1}(t_1, t_2)$ is closed (resp.\ $F_{\sigma}$-set) in
$X$, i.e., $h$ is a contra-continuous (resp.\ Baire-one) function
on $X$.
\end{proof}

The above proof used the technique of Theorem 1 of [4].

\section{Applications}

\begin{dfn}
A real-valued function $f$ defined on a space $X$ is called {\it upper 
semi-contra-continuous} (resp.\ {\it lower semi-contra-continuous}) if 
$f^{-1}(-\infty, t)$ (resp.\ $f^{-1}(t, +\infty)$) is closed for any real 
number $t$.
\end{dfn}

\begin{dfn}
A real-valued function $f$ defined on a space $X$ is called {\it upper 
semi-Baire-one} (resp.\ {\it lower semi-Baire-one}) if $f^{-1}(-\infty, t)$ 
(resp.\ $f^{-1}(t, +\infty)$) is $F_{\sigma}$-set for any real number $t$.
\end{dfn}

The abbreviations $\mathit{usc}, \mathit{lsc}, \mathit{uscc}, 
\mathit{lscc}, \mathit{us}B_1$ and $\mathit{ls}B_1$ are 
used for upper semicontinuous, lower semicontinuous, upper 
semi-contra-continuous, lower semi-contra-continuous, upper 
semi-Baire-one, and lower semi-Baire-one, respectively.

\begin{cor} %\noindent {\bf Corollary 1.} 
Let $g$ and $f$ be real-valued functions defined on a space $X$, that 
$\Lambda$-sets in $X$ are open, such that $f$ is $\mathit{lscc}$, $g$ is 
$\mathit{uscc}$, and $g\leq f$. If $X$ is a extremally disconnected 
space, then there exists a contra-continuous function $h$ defined on $X$ 
such that $g\leq h\leq f$.
\end{cor}

\begin{proof}
Let $g$ be $\mathit{uscc}$, let $f$ be $\mathit{lscc}$, and $g\leq f$.
If a binary relation $\rhobin$ is defined by $A \rhobin B$ in case
$A^{\Lambda}\subseteq B^{V}$, and if $X$ is a extremally
disconnected space then $\rhobin$ is a strong binary relation in the
power set of $X$. For each $t$ in $\mathbb{Q}$, let $A(f,t)$ and
$A(g,t)$ be any lower indefinite cut sets for $f$ and $g$
respectively. If $t_1$ and $t_2$ are any elements of $\mathbb{Q}$
with $t_1<t_2$, then
$$
A(f,t_1)
\subseteq \{x\in X: f(x)\leq t_1\} 
\subseteq \{x\in X: g(x)<t_2\} 
\subseteq A(g, t_2);
$$
since $\{x\in X: f(x)\leq t_1\}$ is open and since $\{x\in X: g(x)<t_2\}$ 
is closed, it follows that $A(f,t_1)^{\Lambda}\subseteq A(g,t_2)^{V}$. 
Hence $t_1<t_2$ implies that $A(f,t_1) \rhobin A(g, t_2)$. 
The proof follows from Theorem 1.
\end{proof}

\begin{cor} %\noindent{\bf Corollary 2.} 
Let $g$ and $f$ be real-valued functions defined on a space $X$, that 
$\Lambda$-sets in $X$ are open, such that $f$ is $\mathit{lscc}$, $g$ is 
$\mathit{uscc}$, and $f\leq g$. If $X$ is a normal space, then there 
exists a contra-continuous $h$ defined on $X$ such that $f\leq h\leq g$.
\end{cor}

\begin{proof}
Let $f$ be $\mathit{lscc}$, $g$ be $\mathit{uscc}$, and $f\leq g$.
A binary relation $\rhobin$ is defined by $A \rhobin B$ in case 
$A^{\Lambda}\subseteq F\subseteq F^{\Lambda}\subseteq B^{V}$ for some 
closed set $F$ in $X$. If $X$ is normal, then $\rhobin$ is a strong binary 
relation in the power set of $X$. If $t_1$ and $t_2$ are any elements of 
$\mathbb{Q}$ with $t_1<t_2$, then
$$
A(g, t_1)
= \{x\in X: g(x)<t_1\}
\subseteq \{x\in X: f(x)\leq t_2\}
=  A(f, t_2);
$$
since $\{x\in X: g(x)<t_1\}$ is closed and since 
$\{x\in X: f(x)\leq t_2\}$ is open and $X$ is normal, it follows that 
$A(g, t_1) \rhobin A(f, t_2)$. 
The proof follows from Theorem 1.
\end{proof}

\begin{cor} % {\bf Corollary 3.} 
Let $g$ and $f$ be real-valued functions defined on a space $X$, that 
$\Lambda$-sets in $X$ are $G_{\delta}$-sets, such that $f$ is 
$\mathit{ls}B_1$, $g$ is $\mathit{us}B_1$, and $g\leq f$. If, for each 
pair of disjoint $G_{\delta}$-sets $G_0, G_1$, there are two 
$F_{\sigma}$-sets $F_0$ and $F_1$ such that $G_0\subseteq F_0$, 
$G_1\subseteq F_1$ and $F_0\cap F_1=\varnothing$, then there exists a 
Baire-one function $h$ defined on $X$ such that $g\leq h\leq f$.
\end{cor}

\begin{proof}
Let $f$ be $\mathit{ls}B_1$, $g$ be $\mathit{us}B_1$, and $g\leq f$. If a 
binary relation $\rhobin$ is defined by $A \rhobin B$ in case 
$A^{\Lambda}\subseteq B^V$, then by hypothesis $\rhobin$ is a strong
binary relation in the power set of $X$. If $t_1$ and $t_2$ are
any elements of $\mathbb{Q}$ with $t_1<t_2$, then
$$
A(f, t_1)\subseteq \{x\in X:  f(x)\leq t_1\}\subseteq \{x\in X: g(x)<
t_2\}\subseteq A(g, t_2);
$$
since $\{x\in X: f(x)\leq t_1\}$ is $G_{\delta}$-set and since
$\{x\in X: g(x)< t_2\}$ is $F_{\sigma}$-set, it follows that 
$A(f,t_1)^{\Lambda}\subseteq A(g, t_2)^{V}$. Hence $t_1<t_2$ implies
that $A(f,t_1) \rhobin A(g, t_2)$. 
The proof follows from Theorem 1.
\end{proof}

\begin{cor} %{\bf Corollary 4.} 
Let $g$ and $f$ be real-valued functions defined on a space $X$, that 
$\Lambda$-sets in $X$ are $G_{\delta}$-sets, such that $f$ is 
$\mathit{ls}B_1$, $g$ is $\mathit{us}B_1$, and $f\leq g$. If, for each 
pair of disjoint $F_{\sigma}$-sets $F_0, F_1$, there are two 
$G_{\delta}$-sets $G_0$ and $G_1$ such that 
$F_0\subseteq G_0, F_1\subseteq G_1$ and $G_0\cap G_1=\varnothing$, 
then there exists a Baire-one function $h$ defined on $X$ such that
$f\leq h\leq g$.
\end{cor}

\begin{proof}
Let $f$ be $\mathit{ls}B_1$, $g$ be $\mathit{us}B_1$, and $f\leq g$.
If a binary relation $\rhobin$ is defined by $A \rhobin B$ in case 
$A^{\Lambda}\subseteq F\subseteq F^{\Lambda}\subseteq B^V$ for some 
$F_{\sigma}$-set $F$ in $X$, then by hypothesis $\rhobin$ is a strong 
binary 
relation in the power set of $X$. If $t_1$ and $t_2$ are any elements of 
$\mathbb{Q}$ with $t_1<t_2$, then
$$
A(g, t_1)
= \{x\in X : g(x)<t_1\} 
\subseteq \{x\in X: f(x)\leq t_2\} 
= A(f, t_2);
$$
since $\{x\in X: g(x)<t_1\}$ is a $F_{\sigma}$-set and since 
$\{x\in X: f(x)\leq t_2\}$ is a $G_{\delta}$-set, by hypothesis it follows 
that $A(g, t_1) \rhobin A(f, t_2)$. 
The proof follows from Theorem 1.
\end{proof}

\begin
{rem}
[{[4], [5]}]
If $g$ and $f$ be real-valued functions defined on a normal space $X$ such 
that $f$ is $\mathit{lsc}$, $g$ is $\mathit{usc}$, and $g\leq f$, then 
there exists a  continuous function $h$ defined on $X$ such that $g\leq  
h\leq f$.
\end{rem}

\begin
{rem}
[{[8]}]
If $g$ and $f$ be real-valued functions defined on a extremally   
disconnected space $X$ such that $f$ is $\mathit{lsc}$, $g$ is   
$\mathit{usc}$, and $f\leq g$, then exists a continuous function $h$  
defined on $X$ such that $f\leq h\leq g$.
\end{rem}

We conclude with the observation that in each of the preceding
corollaries, the separation property for $X$, which is used in order to
verify that $\rhobin$ is a strong binary relation, is also a necessary
condition for the stated insertion property.

\end{document}